\documentclass{article}
\usepackage{epsfig}
\usepackage{amssymb} 
\usepackage{amsthm}
\usepackage{amsmath}

\usepackage{color}
\newcommand{\cb}{\color{black}}

\newcommand{\cn}{\color{black}}

\newtheorem{theorem}{Theorem}
\newtheorem{lemma}{Lemma}

\newcommand{\Bilder}{y}  

\newtheorem{remark}{Remark}

\newcounter{eqsavea}
\newcounter{eqsaveb}

\renewcommand{\rho}{\varrho}
\renewcommand{\phi}{\varphi}
\renewcommand{\theta}{\vartheta}
\newcommand{\nn}{\nonumber}
\newcommand{\pat}{\partial_t}
\newcommand{\pad}{\partial_d}
\newcommand{\pae}{\partial_\varepsilon}

\newcommand{\ve}{\varepsilon}
\newcommand{\half}{\frac12}
\newcommand{\cA}{${\cal A}$}
\newcommand{\cC}{${\cal C}$}
\newcommand{\cG}{{\cal G}}
\newcommand{\cW}{{\cal W}}
\newcommand{\sigmaext}{\sigma_\mathrm{ext}}

\newcommand{\Id}{\mathrm{Id}}
\newcommand{\Wext}{W_\mathrm{ext}}
\newcommand{\whW}{\widehat W}

\newcommand{\ovs}{\overline{\sigma}}
\newcommand{\Div}{\mathrm{div}}
\newcommand{\Det}{\mathrm{det}}

\newcommand{\uve}{\underline{\varepsilon}}
\newcommand{\bn}{\mathbf{n}}
\newcommand{\bu}{\mathbf{u}}
\newcommand{\td}{\tilde{d}}
\newcommand{\tu}{\tilde{\mathbf{u}}}
\newcommand{\te}{\tilde{\varepsilon}}
\newcommand{\teta}{\tilde{\eta}}
\newcommand{\dx}{\,\mathrm{d}x}
\newcommand{\dt}{\,\mathrm{d}t}
\newcommand{\dS}{\,\mathrm{d}S}
\newcommand{\io}{\int\limits_\Omega}
\newcommand{\ioT}{\int\limits_{\Omega_{\cal T}}}
\newcommand{\ioto}{\int\limits_{\Omega_{t_0}}}
\newcommand{\OT}{\Omega_{\cal T}}
\newcommand{\N}{\mathbb{N}} 
\newcommand{\R}{\mathbb{R}} 

\newcommand{\slintO}{\int_\Omega\hspace{-1.25em}-\hspace{0.25em}}
\newcommand{\slintOp}{\int_{\Omega'}\hspace{-1.5em}-\hspace{0.5em}}

\begin{document}

\title{On the Allen-Cahn/Cahn-Hilliard system with a 
geometrically linear elastic energy}
\author{Thomas Blesgen\footnote{Max Planck Institute for Mathematics in
  the Sciences, Inselstra\ss e 22, D-04103 Leip\-zig, Germany,
  email: {\tt blesgen@mis.mpg.de}}~ and
  Anja Schl{\"o}merkemper\footnote{University of W\"urzburg,
  Institute for Mathematics,
  Emil-Fischer-Stra\ss e 40, D-97074 W\"urzburg, Germany,
  email: {\tt anja.schloemerkemper@mathematik.uni-wuerzburg.de}}~}
\date{31 August 2012}

\maketitle

\begin{abstract}
We present an extension of the Allen-Cahn/Cahn-Hilliard system which
incorporates a geometrically linear ansatz for the elastic energy of the
precipitates. The model contains both the elastic Allen-Cahn system and the
elastic Cahn-Hilliard system as special cases and accounts for the
microstructures on the microscopic scale. We prove the existence of weak
solutions to the new model for a general class of energy functionals. We then
give several examples of functionals that belong to this class. This includes
the energy of geometrically linear elastic materials for $D<3$. Moreover we
show this for $D=3$ in the setting of scalar-valued deformations, which
corresponds to the case of anti-plane shear. All this is based on explicit
formulas for relaxed energy functionals newly derived in this article for $D=1$
and $D=3$. In these cases we can also prove uniqueness of the weak solutions.
\end{abstract}

\section{Introduction}
\label{secintro}
In this article we study the Allen-Cahn/Cahn-Hilliard model (AC-CH model for
short), which combines and extends two famous diffuse interface models:
the Allen-Cahn equation and the Cahn-Hilliard equation. Both
have been applied successfully to model segregation, precipitation and
phase change phenomena in alloys and liquid mixtures in materials science,
geology, physics, and biology, among others.

The AC-CH system was first introduced
in \cite{CNC94}, and its mathematical properties have been studied extensively,
see, e.g., \cite{BGN07,NC00} and references therein. We study here an extension
of this model with a particular ansatz for the elastic energy.

The Allen-Cahn equation was first introduced without elasticity in \cite{AC79};
and the Cahn-Hilliard equation was first introduced without elasticity in
\cite{CH}. The Allen-Cahn equation is a second order partial
differential equation of Ginzburg-Landau type for an {\it unconserved}
order-parameter \cb $b$ \cn and thus can be used to model segregation and
precipitation in solids, or other more general situations where a reordering of
the crystal lattice occurs.  Conversely, the Cahn-Hilliard equation is a 
fourth-order partial differential equation for a {\it conserved} order
parameter \cb $a$ \cn striving to model phase change phenomena where
one physical quantity like the volume of the phases is preserved. \cb
Throughout the paper we set $d=a+b$ as this plays a special role. \cn 

The Allen-Cahn system with
linear elasticity was studied before in \cite{BW05}, the Cahn-Hilliard
system with linear elasticity in \cite{Harald1}, \cite{Onuki} and \cite{CL}.
An extension of the Cahn-Hilliard system with geometrically linear elasticity
valid for single crystals was recently found in \cite{BC11} for $D\le2$.

Except for the latter work, elastic effects due to small scale microstructures
within the phases have been neglected. Here we treat the combined AC-CH model
with elasticity in $D\le3$ dimensions. We state the new extended model in
\eqref{ACH1}--\eqref{ACH3} in Section~\ref{secACCH}. This model provides a
basis for
the generalization of further existing isothermal diffuse interface models.

Our definition of the extended model does not require any special assumptions
(except regarding the regularity) on the stored energy functional $\cW$.
However, for the proof of existence of weak solutions,
we require the Assumption~(\cA) phrased in Section~\ref{secexist} below.

We show that Assumption~(\cA) is satisfied in the following cases: 
\begin{itemize}
\item[(i)] Materials which follow the linear theory of elasticity
developed by Eshelby \cite{Eshelby} in the context of elastic inclusions and
inhomogeneities (see a remark after Assumption~(\cA));
\item[(ii)] materials in $D\le2$ which are well-described by a geometrically
linear theory of elasticity (Theorem~\ref{thmcor1}, for $D=2$, cf.\ also
\cite{BS});
\item[(iii)] materials in $D=3$ which are well-described by a geometrically
linear theory of elasticity, where the deformations are assumed to be
scalar-valued functions. This corresponds to the situation of anti-plane shear,
see below for details and cf.\ Theorem~\ref{thmcor1}.
\end{itemize}
The geometrically linear theory of elasticity allows for fine microstructures
within the phases of the elastic materials taken into account. This becomes
apparent in the explicit formulas for the relaxed energy functionals, which
are, together with formulas for their derivatives, the basis of our proof of
Theorem~\ref{thmcor1}. In Section~\ref{secexpl} we recall an explicit formula
in $D=2$ derived by Chenchiah and Bhattacharya \cite{CB2008}. Moreover we prove
explicit formulas for $D=1$, and for $D=3$ in the setting of scalar-valued
deformations. (In \cite{CB2008} partial results for $D=3$ in the vectorial,
i.e. the non-scalar, setting are shown.)  
 
The coupling to elasticity changes significantly the morphology of the
precipitates and the coarsening patterns, see, e.g., the classification in
\cite{FPL99}. For an intuitive picture of the microstructures taken into
account we refer the reader to Figures~\ref{fig:intuitivemacro} and
\ref{fig:intuitivemicro}. 
\begin{figure}
\begin{center} \input macro.pstex_t
\caption{Macroscopic phases of a segregated material with phases $d$ being
almost one (dark gray) and almost zero (light gray). In
Figure~\ref{fig:intuitivemicro} blow-ups of the circular regions
are shown.}
\label{fig:intuitivemacro}
\end{center}
\end{figure}
\begin{figure}
\begin{center} \input micro2.pstex_t \qquad 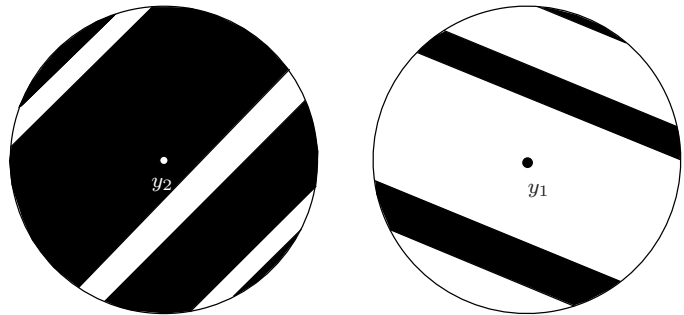
\caption{Microstructures with $d$ being almost one (left) and almost zero
(right). The black part
corresponds to $\tilde d =1$, while $\tilde d =0$ in the white
regions. The average of $\tilde d$ equals $d$.}
\label{fig:intuitivemicro}
\end{center}
\end{figure}
Figure~\ref{fig:intuitivemacro} displays exemplary macroscopic phases of a
segregated material with phases $d$ being almost one (dark gray) and almost zero
(light gray). In Figure~\ref{fig:intuitivemicro}
we show the corresponding microscopic length scale featuring
microstructures; the left is a blow-up
of a small region within the phase with $d$ being almost one and the right is a
blow-up of a very small region within the phase with $d$ being almost zero.
On this microscopic scale where we assume that it is sufficient to treat the
elastic energies of single crystals, fine microstructures in the form of
laminates occur, see, e.g., \cite{Salje, Gottstein}.
 
The methods developed here apply generally to any established phase change and
segregation model provided the temperature is conserved. (For non-isothermal
settings, the validity of the second law of thermodynamics requires additional
corrections which are not studied here.) For further discussions of our model
and the analytical results we refer to the conclusions in
Section~\ref{secconclude}. 

\section{The AC-CH model and extensions}
\label{secACCH}
Throughout this paper, let $\Omega\subset\R^D$ for $D\ge1$ be a bounded domain
with Lipschitz boundary which serves as an (unstressed) reference configuration.
For a stop time ${\cal T}>0$, let $\OT:=\Omega\times(0,{\cal T})$ denote the
space-time cylinder. To the Allen-Cahn/Cahn-Hilliard system, first derived in
\cite{CNC94}, we add elasticity, possibly respecting the lamination
microstructure, by introducing the system 
\setcounter{eqsavea}{\value{equation}} 
\begin{eqnarray}
\label{ACH1}
\pat a &=& \lambda\,\Div\Big(M(a,b)\nabla\frac{\partial F}{\partial a}\Big),\\
\label{ACH2}
\pat b &=& -M(a,b)\,\frac{\partial F}{\partial b},\\
\label{ACH3}
\mathbf{0} &=& \Div\left(\partial_{\ve}\cW(a+b,\ve(\bu))\right),
\end{eqnarray}
where the function $a:\OT\to\R_0^+$ is a
{\it conserved} order parameter, typically a concentration, $b:\OT\to\R_0^+$ is
an {\it unconserved} order-parameter, specifying the reordering of the
underlying lattice, $M(a,b)$ denotes the positive semi-definite mobility
tensor, $\lambda>0$ is a small constant determining the interfacial thickness.
The choice of $\cW$ determines whether lamination microstructure occurs in the
system, see \eqref{Wmin}.

By $\bu:\Omega\to\R^D$ we describe the displacement field, such that a material
point $x$ in the undeformed body $\Omega$ is at $x'=x+\bu(x)$ after the
deformation. Then the (linearized) strain tensor is defined by
\begin{equation}
\label{vedef}
\ve(\bu):=\half\left(\nabla\bu+\nabla\bu^t\right),
\end{equation}
where $A^t$ denotes the transpose of a matrix $A\in \R^{D\times D}$.
As usual, $\cdot$ stands for the inner product in $\R^D$, that is
$\bu\cdot\mathbf{v}=\sum_{i=1}^D u_iv_i$, and for $A,\,B\in\R^{D\times D}$ we
denote the inner product in $\R^{D\times D}$ by 
\[ A\!:\!B:=\mathrm{tr}(A^tB)=\sum_{i,j=1}^DA_{ij}B_{ij}. \]
Moreover, $|A|:=\sqrt{A\!:\!A}$ for
$A\in\R^{D\times D}$ is the Frobenius norm. We denote the symmetric matrices
in $\R^{D\times D}$ by $\R^{D\times D}_{\operatorname{sym}}$.

The functional
$\cW(a+b,\ve(\bu))$ represents the stored elastic energy density. We will choose
it either according to the linear ansatz by Eshelby, $\cW=W_\mathrm{lin}$, or
as a single crystal composite lamination energy, $\cW=\whW$. 
The definition of $\whW$ is given in the coming section \ref{ssectelast}.

The linear theory by Eshelby \cite{Eshelby} developed in the
context of elastic inclusions and inhomogeneities, can be summarized
in the following ansatz for the elastic energy
\begin{equation}
\label{Wlin}
W_{\mathrm{lin}}(d,\ve):=\half(\ve-\uve(d)):C(d)(\ve-\uve(d))
\end{equation}
for all $\ve\in\R_{\operatorname{sym}}^{D\times D}$, $d:=a+b$, and
$\uve(d):=d\,\uve$ with a constant
$\uve\in\R_{\operatorname{sym}}^{D\times D}$. 

By $C(d)$ we denote the symmetric, positive definite and
concentration-dependent elasticity tensor of the system that maps
symmetric tensors in $\R^{D\times D}$ to themselves.

The system (\ref{ACH1})--(\ref{ACH3}) is completed with the definition of
the free energy
\begin{equation}
\label{Fdef}
F(a,b,\bu):=\io\psi(a,b)+\frac{\lambda}{2}\left(|\nabla a|^2+|\nabla b|^2\right)
+\cW(a+b,\ve(\bu))+\Wext(\ve(\bu))\dx,
\end{equation}
see \cite{CNC94}, where $\psi(a,b)$ is the free energy density assumed to be
\begin{eqnarray}
\label{psidef}
\psi(a,b) &:=& \frac{\theta}{2}\Big(g(a+b)+g(a-b)\Big)+\kappa_1a(1-a)
-\kappa_2 b^2,\\
g(s) &:=& s\ln s+(1-s)\ln(1-s)\nn
\end{eqnarray}
for scalars $\kappa_1,\,\kappa_2>0$.
The term $\half[g(a+b)+g(a-b)]$ in (\ref{psidef}) defines the entropic part of
the free energy, given in the canonical Bernoulli form for perfect mixing, and
$\theta>0$ is the constant temperature.

The functional $\Wext(\ve)$ in (\ref{Fdef}) represents energy effects due to
applied forces. In the absence of body forces, the work necessary to transform
the undeformed body $\Omega$ into a state with displacement $\bu$ is then
\[ -\int_{\partial\Omega}\bu\cdot\sigmaext\bn\dS=-\io\nabla\bu:\sigmaext\dx
=-\io\ve(\bu):\sigmaext\dx, \]
where we use that the applied stress $\sigmaext$ is constant and symmetric.
Consequently,
\begin{equation}
\label{Wstar} \Wext(\ve)=-\ve:\sigmaext
\end{equation}
is the energy density of the applied outer forces. The system
(\ref{ACH1})--(\ref{ACH3})
has to be solved in $\OT$ subject to the initial conditions
\[ a(t=0)=a_0,\; b(t=0)=b_0\quad\mbox{in }\Omega \]
for given functions $a_0,\,b_0: \Omega\to\R$ subject to the Neumann boundary
conditions for $a$, the no-flux boundary conditions, and the equilibrium
condition for applied forces
\begin{equation}
\label{BC}
\nabla a\cdot\bn \cb = \nabla b \cdot \bn \cn =0,\quad J(a,b,\bu)\cb \cdot \bn
\cn=0,\quad\sigma\cdot\bn=
\sigmaext\cdot\bn\qquad\mbox{on }\partial\Omega,\, t>0.
\end{equation}
Here, $\sigma:=\pae\cW(a+b,\ve(\bu))$ defines the stress.

In (\ref{BC}), $\bn$ is the unit outer normal to $\partial\Omega$.
For simplicity, body forces are neglected and it is assumed that the boundary
tractions are dead loads given by a constant symmetric tensor $\sigmaext$.
By $J$ we denote the mass flux, given by
\[ J(a,b,\bu):=-M(a,b)\nabla\mu=-M(a,b)\nabla
\frac{\partial F}{\partial a}(a,b,\bu), \]
with $\mu:=\frac{\partial F}{\partial a}$ the chemical potential.

The valid parameter range of $a$ and $b$ is, see Theorem~\ref{theo1},
\begin{equation}
\label{abrange} 0\le a+b\le1,\qquad 0\le a-b\le 1.
\end{equation}
The inequalities are strict unless $(a,b)=(0,0)$ or $(a,b)=(1,0)$.

\cb For $a_0=b_0\equiv0$ in $\Omega$, we obtain the pathological solution
$a=b\equiv0$ in $\OT$. For $a_0\equiv1$, $b_0\equiv0$ in $\Omega$ we
obtain a Cahn-Hilliard equation in $b$ and a pathological equation for $a$. \cn

The system (\ref{ACH1})--(\ref{ACH3}) includes as special case the
{\it elastic Cahn-Hilliard system} (setting $b\equiv0$, \cite{Harald1}) 
\begin{eqnarray*}
\pat a &=& \lambda\,\Div\Big(M(a)\nabla\frac{\partial F}{\partial a}\Big),\\
\mathbf{0} &=& \Div\left(\partial_{\ve}\cW(a,\ve(\bu))\right)
\end{eqnarray*}
with 
\begin{align*}
F(a,\bu) &:= \io\psi(a)+\frac{\lambda}{2}|\nabla a|^2
+\cW(a,\ve(\bu))+\Wext(\ve(\bu))\dx,\\
\psi(a) &:= \frac{\theta}{2}a\ln a+(1-a)\ln(1-a)+\kappa_1a(1-a)
\end{align*}
and the boundary and initial conditions correspondingly to above. Moreover the
system (\ref{ACH1})--(\ref{ACH3}) includes as special case the
{\it elastic Allen-Cahn equations} (setting $a\equiv\half$, \cite{BW05}) which
for rescaled $b$ with $0<b<1$ read
\begin{eqnarray*}
\pat b &=& -M(b)\,\frac{\partial F}{\partial b},\\
\mathbf{0} &=& \Div\left(\partial_{\ve}\cW(b,\ve(\bu))\right),
\end{eqnarray*}
with
\begin{align*}
F(b,\bu) &:= \io\psi(b)+\frac{\lambda}{2}|\nabla b|^2+\cW(b,\ve(\bu))
+\Wext(\ve(\bu))\dx,\\
\psi(b) &:= \frac{\theta}{2}\Big(g(b)+g(1-b)\Big)+\frac{\kappa_1}4-\kappa_2b^2
\end{align*}
and $g$ as well as the boundary and initial conditions are as above.
\pagebreak

The system (\ref{ACH1})--(\ref{ACH3}) is exemplary for an isothermal
model that exhibits simultaneous ordering and phase transitions.
Equation~(\ref{ACH1}) is a diffusion law for $a$ governed by the flux $J$
and states the conservation of mass in $\Omega$. Equation~(\ref{ACH2}) is a
simple gradient flow in the descent direction $-\frac{\partial F}{\partial b}$.
Equation~(\ref{ACH3}) is a consequence of Newton's second law under the
additional assumption that the acceleration $\partial_{tt}\bu$ originally
appearing on the left hand side can be neglected (this can be proved formally
by a scaling argument and formally matched asymptotics). The vector equation
(\ref{ACH3}) serves to determine the unknown displacement $\bu$. 

\begin{remark}
\label{rem1}
The equations~(\ref{ACH1})--(\ref{ACH3}) can be generalized to vector-valued
mappings $a$, $b$. This allows to study situations with more than two phases
present. To fix ideas and for the sake of a clear presentation, we restrict
ourselves throughout this paper to scalars $a$ and $b$.
\end{remark}

\begin{remark}
\label{rem2}
The equations~(\ref{ACH1})--(\ref{ACH3}) with boundary conditions (\ref{BC})
and \eqref{Fdef} and \eqref{psidef} comply with the second law of
thermodynamics, which in case of isothermal conditions reads for a closed system
\[ \pat F(a(t),b(t),\bu(t))\le0. \]
This inequality can be verified by direct inspection similar to the
calculations in \cite{BC11}.
\end{remark}

In the proof of Theorem~\ref{theo2} we apply the following explicit
formulation of (\ref{ACH1})--(\ref{ACH3}) with constant mobility $M\equiv1$
\setcounter{eqsaveb}{\value{equation}}
\renewcommand{\theequation}{\arabic{section}.\arabic{equation}'}
\setcounter{equation}{\value{eqsavea}}
\begin{eqnarray}
\label{ACH1e}
\pat a &=& \lambda\triangle\!\left[\frac{\theta}{2}\big(g'(a\!+\!b)+g'(a\!-\!b)
\big)+\kappa_1(1\!-\!2a)+\frac{\partial\cW}{\partial d}(a\!+\!b,\ve(\bu))
-\triangle a\!\right]\!\!,\\
\label{ACH2e} \pat b &=& \lambda\triangle b+\frac{\theta}{2}\left[g'(a-b)
-g'(a+b)\right]+2\kappa_2b-\frac{\partial\cW}{\partial d}(a+b,\ve(\bu)),\\
\label{ACH3e} \mathbf{0} &=& \Div\left(\pae\cW(a+b,\ve(\bu))\right).
\end{eqnarray}
\setcounter{equation}{\value{eqsaveb}}
\renewcommand{\theequation}{\arabic{section}.\arabic{equation}}

The subsequent section is devoted to the geometrically linear elastic energy
for single crystals. This is a prerequisite to the discussion of existence and
uniqueness results for the extensions of the Allen-Cahn/Cahn-Hilliard
(AC-CH) model studied in Section~\ref{secexist}.

\subsection{The geometrically linear theory of elasticity in single crystals}
\label{ssectelast}
Our main objective in this subsection is to study a geometrically linear theory
of elasticity in the context of isothermal phase transitions.
For systematic reasons, we first recall the linear ansatz dating back to
Eshelby, \cite{Eshelby}. As a byproduct of the existence theory proved in
Section~\ref{secexist} we obtain a new existence
result for the AC-CH equations with linear elasticity. Subsequently, we
introduce the geometrically linear elasticity theory that takes the laminates
of the material into account. 

\vspace*{2mm}
In the following we assume that two phases are present in the considered
material which may form microstructures as displayed, e.g., in
Figures~\ref{fig1} and \ref{fig2}. We refer to the energy $W_i$, $i=1,2$ of
each of the phases as {\it microscopic energy}, cf.\ \eqref{Widef}, and to the
energy $\whW(d,\ve(\bu))$ in \eqref{Wmin}, which reflects the effective
behavior of the system with microstructures, as the {\it macroscopic energy}.

To determine the energy $\whW(d,\ve(\bu))$ in the geometrically linear theory we
need to solve a local minimization problem, Eqn.~(\ref{Wmin}) below, which we
shall outline now.

We assume that the volumes occupied by each of the two phases in $\Omega$ are
measurable sets. In particular, if $\td_1\equiv\td$, $\td_2\equiv1-\td$
characterize the two phases on the microscale, we have
$\td_i\in BV(\Omega;\{0,1\})$ and $\td_1+\td_2=1$ a.e. in $\Omega$. The symbol
$BV$ denotes the space of functions of bounded variation, see, e.g.,
\cite{AFP,Ziemer}. By
\begin{equation}
\label{angledef}
\langle\tilde{\rho}\rangle:=\slintO\tilde{\rho}(x)\dx:=\frac{1}{|\Omega|}\io
\tilde{\rho}(x)\dx
\end{equation}
we denote the average of a function $\tilde{\rho}$ in $\Omega$, where
$|E|$ is the $D$-dimensional Lebesgue measure of a set $E$.

Let $\ve_i^T\in\R_{\operatorname{sym}}^{D\times D}$, $i=1,2$, be
the stress-free strain (or eigenstrain) of the $i$-th phase relative to the
chosen reference configuration and $\alpha_i$ be its positive definite
elasticity tensor. Then the elastic energy density of phase $i$ subject to a
strain $\te$ is given by
\begin{equation}
\label{Widef}
W_i(\te):=\half\alpha_i\left(\te-\ve_i^T\right):\left(\te-\ve_i^T\right)+w_i
\end{equation}
for $w_i\ge0$. In (\ref{Widef}) we assumed for simplicity that $\alpha_i$,
$\ve_i^T$ and $w_i$ are constants, independent of $x\in\Omega$ and the order
parameter $d$.

Under the assumption that the elastic energy adapts infinitely fast and that
the surface energy between laminates of the microstructure can be neglected,
the effective elastic energy is, \cite{CB2008},
\begin{equation}
\label{Wmin} \whW(d,\ve)(y):=\!\inf_{\substack{\langle\td\rangle=d(y)\\
\td\in\{0,1\}}}\;\inf_{\tu_{|\partial\Omega}=\ve(y)x}\slintO W(\td,\ve(\tu))\dx,
\quad d\in[0,1],
\end{equation}
where we used
\begin{equation}
\label{WWdef}
W(\td,\te):=\td W_1(\te)+(1-\td)W_2(\te),\quad\td\in\{0,1\}.
\end{equation}
The definition (\ref{Wmin}) requires further clarification. Firstly,
$\te=\ve(\tu):=\half(\nabla\tu+\nabla\tu^t)$, and instead of integrating
over $\Omega$, one may integrate over $B_r(y)$, the open ball of radius $r$
around $y\in\Omega$, where the mean $\langle\cdot\rangle$ is now taken over
$B_r(y)$. By homotopy arguments or by results in \cite{Dac}, any $r>0$ yields
the same value of $\whW(d,\ve)(y)$, as long as $B_r(y)\subset\Omega$. Taking
the union of such balls then leads to (\ref{Wmin}). Secondly,
the infimum over $\td$ is the result of homogenization subject to the
constraint that the volume fraction of the selected phase is preset by $d(y)$,
see \cite[Chapter~10]{Donato}.
This infimum is taken over functions $\td\in BV(\Omega;\,\{0,1\})$ as explained
above ensuring $\whW\ge0$. This is why (\ref{Wmin}) is only meaningful for
$d\in[0,1]$.

The second infimum is taken over functions $\tu\in H^1(\Omega;\,\R^D)$ where
the condition
$\tu_{|\partial\Omega}=\ve(y)x$ has to be read as $\tu(x)=\ve(y)x$ for a.e.
$x\in\partial\Omega$. This originates from the requirement that the functional
$\whW$ thus defined must be quasi-convex,
see \cite{Dac}, and is the result of
relaxation theory, \cite{Dac}, \cite{Kohn}, as follows. If for prescribed
$d=a+b$ the microscopic elastic energy density is denoted by $W_d(\te)$, then
\begin{equation}
\label{Wpre}
\whW_d(\ve):=\inf_{\tu_{|\partial\Omega}=\ve x}\;\slintO W_d(\ve(\tu))\dx
\end{equation}
is the elastic energy density of the material with macroscopic strain $\ve$
after microstructure has formed. For $D=2$, explicit analytic formulas for
$\whW$ are known, see \cite{CB2008} as well as formulas of the partial
derivatives of $\whW$. These will be recalled in Section~\ref{secexpl}.

\section{Existence and uniqueness results for the AC-CH system} \label{secexist}
The existence of solutions to the Allen-Cahn/Cahn-Hilliard equation without 
elasticity was studied in \cite{BCN94} with the help of a semigroup calculus.
Existence and uniqueness of weak solutions to the Cahn-Hilliard equation with
linear elasticity is proved in \cite{Harald1}, with geometrically linear
elasticity in \cite{BC11}. Existence and uniqueness of weak solutions to the
Allen-Cahn equation with linear elasticity is shown in \cite{BW05}.

Subsequently we provide existence and uniqueness results for
(\ref{ACH1})--(\ref{ACH3}), where $\cW$ is an elastic energy density
satisfying the following assumption (\cA).

\vspace*{2mm}
(\cA) The elastic energy density
$\cW\in C^1(\R\times\R^{D\times D}_\mathrm{sym};\,\R)$ satisfies the conditions

(\cA1) $\pae\cW(d,\cdot)$ is strongly monotone uniformly in $d$, i.e., there
exists a constant\\
\phantom{(\cA1) }$c_1>0$ such that for all
$\ve_1,\,\ve_2\in\R_{\operatorname{sym}}^{D\times D}$ and all
$d\in\R$
\[ \left(\pae\cW(d,\ve_2)-\pae\cW(d,\ve_1)\right):\left(\ve_2-\ve_1\right)
\ge c_1|\ve_2-\ve_1|^2. \]
(\cA2) There exists a constant $C_1>0$ such that for all $d\in\R$ 
and all $\ve\in\R_{\operatorname{sym}}^{D\times D}$
\begin{eqnarray*}
|\cW(d,\ve)| &\le& C_1(|d|^2+|\ve|^2+1),\\
|\pad\cW(d,\ve)| &\le& C_1(|d|^2+|\ve|^2+1),\\
|\pae\cW(d,\ve)| &\le& C_1(|d|+|\ve|+1).
\end{eqnarray*}
All constants in this article, unless explicitly stated otherwise, may depend
on the material parameters $\alpha_1$, $\alpha_2$, $\ve_1^T$ and $\ve_2^T$, but
are independent of $d$ and $\ve$.
Condition~(\cA1) states that $\cW$ is convex in $\ve$. The problem
becomes non-convex through the dependence on $d$.
One prominent example satisfying assumption (\cA) is the elastic energy 
$W_\mathrm{lin}$ in \eqref{Wlin}. In Section~\ref{secexam} we prove that also
the relaxed energy functional $\whW$ defined in (\ref{Wmin})
satisfies the assumption (\cA).

\cb We require the condition
\begin{equation}
\label{W0bound}
\cW(a_0+b_0,\ve(\bu(x,0)))<\infty,
\end{equation}
where $\bu(\cdot,0)$ is the solution of \cb (\ref{ACH3}) for
$a=a_0$, $b=b_0$.\cn

\vspace*{2mm}
\begin{theorem}[Existence of weak solutions]
\label{theo1}
Let the mobility tensor $M$ be positive definite \cb and continuous for all
$a$, $b$ satisfying (\ref{abrange}), \cn $\cW$ fulfill (\cA),
$\psi$ be given by (\ref{psidef}) and the initial data $(a_0,b_0)$ satisfy
(\ref{abrange}) \cb and (\ref{W0bound}). \cn
Then there exists a \cb weak \cn solution $(a,b,\bu)$ to
(\ref{ACH1})--(\ref{ACH3}) that satisfies

(i) $a,b\in C^{0,\frac14}\left([0,{\cal T}];\,L^2(\Omega)\right)$,\\
(ii) $\pat a \cb\in L^2(0,{\cal T}; H^1(\Omega)^*)$,
$\,\pat b\in L^2(\OT)$,\\
(iii) $\bu\in L^\infty\left(0,{\cal T};\,H^1(\Omega;\,\R^D)\right)$,\\
(iv) The feasible parameter range of $(a,b)$ is given by (\ref{abrange}).
\end{theorem}

\vspace*{2mm}
{\bf Proof}: The statements of the theorem can be proved with the methods
developed in \cite{BW05} for an Allen-Cahn system with linear elasticity \cb
and in \cite{Harald1} for a Cahn-Hilliard system with linear elasticity. \cn
We sketch the main steps.

First we introduce the operator ${\cal M}$ associated to
$w\mapsto-M\triangle w$ as a mapping from $H^1(\Omega)$ to its dual by 
\begin{equation}
\label{Meq}
{\cal M}(w)\eta:=\io M\nabla w\cdot\nabla\eta\dx.
\end{equation}
From the Poincar{\'e} inequality and the Lax-Milgram theorem 
(which can be applied since $M$ is assumed to be positive
definite) we know that ${\cal M}$ is invertible and we denote its inverse by
$\cG$, the Green's function. We have
\[ (M\nabla\cG f,\nabla\eta)_{L^2}=\langle\eta,f\rangle\quad\mbox{for all }
\eta\in H^1(\Omega),\,f\in(H^1(\Omega))'. \]
For $f_1,\,f_2\in(H^1(\Omega))'$, we define the inner product
\[ (f_1,f_2)_M:=(M\nabla\cG f_1,\nabla\cG f_2)_{L^2} \]
with the corresponding norm
\[ \|f\|_M:=\sqrt{(f,f)_M}\quad\mbox{for }f\in(H^1(\Omega))'. \]

For a small discrete step size $h>0$, chosen such that ${\cal T}h^{-1}\in\N$,
for time steps $m\in\N$ with $0<m<{\cal T}h^{-1}$, and given values
$a^{m-1}$, $b^{m-1}\in\R$, we introduce the discrete free
energy functional
\begin{equation}
\label{Fhdef} F^{m,h}(a,b,\bu):=F(a,b,\bu)+\frac{1}{2h}\|a-a^{m-1}\|_M^2
+\frac{1}{2h}\|b-b^{m-1}\|_{L^2}^2,
\end{equation}
where (in case of $m=1$) it holds $a^0=a_0$, $b^0=b_0$, the initial values
of $a$ and $b$. By the direct method in the calculus of variations and
Assumption~(\cA), it is possible to show that
for $h$ sufficiently small, $F^{m,h}$ possesses a minimizer
$(a^m,b^m,\bu^m)\in H^1(\Omega)\times H^1(\Omega)\times H^1(\Omega;\,\R^D)$.
This minimizer solves the fully implicit time discretization of
(\ref{ACH1})--(\ref{ACH3}). Next the discrete solution is extended
affine linearly to $(\overline{a},\overline{b},\overline{\bu})$ by setting for
$t=(\tau m+(1-\tau)(m-1))h$ with suitable $\tau\in[0,1]$
\[ (\overline{a},\overline{b},\overline{\bu})(t):=\tau(a^m,b^m,\bu^m)
+(1-\tau)(a^{m-1},b^{m-1},\bu^{m-1}). \]
The validity of the second law of thermodynamics (cf.\
Remark~\ref{rem2}) \cb together with (\ref{W0bound}) \cn implies that $F$ is
non-increasing in time. \cb In combination with a higher integrability
condition on $\bu$, \cite{Harald1}, \cn this allows to derive uniform estimates
for
$(\overline{a},\overline{b},\overline{\bu})$. Compactness arguments then allow
to pass to the limit $h\searrow0$ and the limit solves
(\ref{ACH1})--(\ref{ACH3}). \qed

\vspace*{2mm}
In general, the uniqueness of solutions to \eqref{ACH1}--\eqref{ACH3}
is open. However, we prove it in a special case for the linear elastic
energy density $\cW=W_\mathrm{lin}$.

\begin{theorem}[Uniqueness of solutions for linear elasticity]
\label{theo2}$\hspace*{-1.3pt}$Let $\cW\!\!=\!\!W_\mathrm{lin}$ be given by
(\ref{Wlin}), the material be homogeneous, i.e. the elasticity tensor $C$ be
independent of $d$, and let $M\equiv1$. Then the solution $(a,b,\bu)$ of
Theorem~\ref{theo1} is unique in the spaces stated there.
\end{theorem}

\vspace*{2mm}
{\bf Proof}: The proof is very similar to the proof of uniqueness for the
Cahn-Hilliard system \cb \cite{Harald1} \cn but we repeat it here because we
later need to modify it.

Fix $t_0\in(0,{\cal T})$.
Let $(a^k,b^k,\bu^k)$, $k=1,2$ be two pairs of solutions to
(\ref{ACH1e})--(\ref{ACH3e}) and (\ref{Wlin}). The differences $a:=a^1-a^2$,
$b:=b^1-b^2$, $\bu:=\bu^1-\bu^2$ with corresponding difference of the
chemical potentials $\mu:=\mu^1-\mu^2:=\frac{\partial F}{\partial a}(a^1,b^1)
-\frac{\partial F}{\partial a}(a^2,b^2)$ solve the weak equations
\begin{align}
\label{test1}
& \ioT\left[-a\pat\xi+\lambda\nabla\mu\cdot\nabla\xi\right]\dx\dt=0,\\
& \ioT\left[\pat b\eta+\lambda\nabla b\cdot\nabla\eta-\uve:C\left(\ve(\bu)
-(a+b)\uve\right)\eta\right]\dx\dt\nn\\
\label{test2}
& \; =\!\ioT\left[\frac{\theta}{2}\left(g'(a^2\!+\!b^2)-g'(a^1\!+\!b^1)
+g'(a^1\!-\!b^1)-g'(a^2\!-\!b^2)\right)\eta+2\kappa_2b\eta\right]\dx\dt,\\
\label{test3}
& \ioto C\left(\ve(\bu)-\uve(a+b)\right):\ve(\bu)\dx\dt=0
\end{align}
for every $\xi,\,\eta\in L^2(0,{\cal T};\,H_0^1(\Omega))\cap L^\infty(\OT)$ with
$\pat\xi,\,\pat\eta\in L^2(\OT)$, $\xi({\cal T})=0$, where in order to get
(\ref{test3}) we plugged in $(\bu^2-\bu^1){\cal X}_{(0,t_0)}$ as a test function
and integrated by parts. As a test function in (\ref{test1}) we pick
\[ \xi(x,t):=\left\{\!\begin{array}{ll}\int_t^{t_0}\mu(x,s)\,ds,
&\quad\mbox{if }t\le t_0,\\
0, &\quad\mbox{if }t>t_0.\end{array}\right. \]
This shows
\begin{equation}
\label{step1} \ioto a\mu+\lambda\nabla(\cG a)\cdot\nabla(\pat\cG a)\dx\dt=0.
\end{equation}
The difference of the chemical potentials fulfils, with the help of
(\ref{Fdef}),
\begin{align*}
\ioT\mu\zeta\dx\dt &= \ioT\Big[\frac{\theta}{2}\Big(g'(a^1+b^1)-g'(a^2+b^2)
+g'(a^2-b^2)-g'(a^1-b^1)\Big)\zeta\\
& \hspace*{35pt} -2\kappa_1a\zeta+\lambda\nabla a\cdot\nabla\zeta
-\uve:C(\ve(\bu)-(a+b)\uve)\zeta\Big]\dx\dt.
\end{align*}
We pick $\zeta:=(a^1-a^2){\cal X}_{(0,t_0)}$. With (\ref{step1}) we obtain
\begin{align}
& \hspace*{-20pt} \frac{\lambda}{2}\|a(t_0)\|_M^2+\ioto\lambda|\nabla a|^2
-a\uve:C(\ve(\bu)-\uve(a+b))\dx\dt\le\nn\\
\label{step2}
& \ioto 2\kappa_1a^2+\frac{\theta}{2}\Big[\big|g'(a^1\!+\!b^1)
-g'(a^2\!+\!b^2)\big|+\big|g'(a^1\!-\!b^1)-g'(a^2\!-\!b^2)\big|\Big]|a|\dx\dt.
\end{align}
In (\ref{test2}) we choose $\eta:=(b^1-b^2){\cal X}_{(0,t_0)}$ as a test
function and add the resulting equation to (\ref{step2}) and use (\ref{test3}).
We end up with
\begin{eqnarray*}
&& \hspace*{-30pt} \frac{\lambda}{2}\|a(t_0)\|_M^2+\half\|b(t_0)\|_{L^2}+\ioto
\Big[\lambda\big(|\nabla a|^2+|\nabla b|^2\big)+\cW(a+b,\ve(\bu))\Big]\dx\dt\\
&& \le\ioto2(\kappa_1|a|^2+\kappa_2|b|^2)\dx\dt\\
&& \quad+\!\!\ioto\!\frac{\theta}{2}\Big[\big|g'(a^1\!+\!b^1)
-g'(a^2\!+\!b^2)\big|+\big|g'(a^1\!-\!b^1)-g'(a^2\!-\!b^2)\big|\Big]
(|a|\!+\!|b|)\dx\dt.
\end{eqnarray*}
From Theorem~\ref{theo1} we know that the terms $g'(a^i\pm b^i)$, $i=1,2$
are finite, and $g'$ is Lipschitz continuous. Applying first Young's
inequality, then Gronwall's inequality, as $t_0\in(0,{\cal T})$ was arbitrary,
we find $a=b=0$ in $\OT$. This finally yields
\[ \ioT\ve(\bu):C\ve(\bu)\dx\dt=0. \]
With Korn's inequality this proves $\bu\equiv0$ in $\OT$.\qed

\section{Explicit formulas for $\whW$}
\label{secexpl}
In many situations like the numerical implementation of the extended models,
the above definition (\ref{Wmin}) of $\whW$ is not practical since it is
indirect and based on a local minimization. For these applications and for
direct later use, we collect here some explicit formulas of the relaxed energy
$\whW$ for $D\le2$ and for the scalar case in $D=3$.

\subsection{The case $D=2$}
\label{subexpl2d}
As shown in \cite{CB2008}, it holds
\begin{equation}
\label{Wdef}
\whW(d,\ve)=d_1W_1(\ve_1^*)+d_2W_2(\ve_2^*)+\beta^*d_1d_2\det(\ve_2^*-\ve_1^*),
\end{equation}
where $\beta^*$, $\ve_1^*$ and $\ve_2^*$ are defined below.
First we need to fix further notations following \cite{CB2008}.

Let $\gamma^*>0$ be given by
\begin{equation}
\label{gsdef}
\gamma^*:=\min\{\gamma_1,\gamma_2\},
\end{equation}
where $\gamma_i$ is the reciprocal of the largest eigenvalue of
$\alpha_i^{-1/2}T\alpha_i^{-1/2}$, $\alpha_i$ is the elastic modulus of
laminate $i$, and the operator
$T:\R^{2\times2}_\mathrm{sym}\to\R^{2\times2}_\mathrm{sym}$ is given by
\[ T\ve=\ve-\operatorname{tr}(\ve)\Id. \]
In \cite{BC11} a recipe is given for the practical computation of $\gamma^*$.
Here, we only remark that if the space groups of the two existing
laminates are cubic, it holds
\[ \gamma^*=\min\{C_{1,11}-C_{1,12},C_{2,11}-C_{2,12},\,2C_{1,44},
2C_{2,44}\}. \]
The first subscript of $C$ denotes here the phase, the other two indices are
the coefficients of the reduced elasticity tensor in Voigt notation,
\cite{Nye}.

As shown in \cite{CB2008}, the scalar $\beta^*\in[0,\gamma^*]$ determines the
amount of translation of the laminates defined by
\begin{equation}
\label{bsdef}
\beta^*=\beta^*(d,\ve):=\left\{
\begin{array}{l@{\quad\mbox{if }}l@{\quad\mbox{(Regime }}l}
0 & \varphi\equiv0 & 0\mbox{)},\\
0 & \varphi(0,d,\ve)>0 & \mbox{I)},\\
\beta_{II} & \varphi(0,d,\ve)\le0\mbox{ and }\varphi(\gamma^*,d,\ve)\ge0 &
\mbox{II)},\\
\gamma^* & \varphi(\gamma^*,d,\ve)<0 & \mbox{III).} \end{array}\right.
\end{equation}
In this definition, $\beta_{II}=\beta_{II}(d,\ve)$ is the unique
solution of $\varphi(\cdot,d,\ve)=0$ with $\varphi$ defined by
\begin{eqnarray}
\label{phidef}
\varphi(\beta,d,\ve) &:=& -\Det(\triangle\ve^*(\beta,d,\ve)) \,\,=\,\,
-\Det\Big[\alpha(\beta,d)^{-1}e(\ve)\Big],\\
\triangle\ve^* &=& \triangle\ve^*(\beta,d,\ve)
\,\,:=\,\,\ve_2^*(\beta,d,\ve)-\ve_1^*(\beta,d,\ve),\nn
\end{eqnarray}
and the yet undefined functions are specified below.

\vspace*{2mm}
The four regimes have the following crystallographic interpretation, which
follows from the construction of the optimal microstructure in the calculation
of $\whW$.

{\bf Regime~0}: The material is homogeneous and the energy does not depend
on the microstructure.
This occurs when $\alpha_2(\ve_2^T-\ve)-\alpha_1(\ve_1^T-\ve)=0$.

{\bf Regime~I}: There exist two optimal rank-I laminates.

This is characterized by
\[ \det\big[(d_2\alpha_1+d_1\alpha_2)^{-1}\big(\alpha_2(\ve_2^T-\ve)
-\alpha_1(\ve_1^T-\ve)\big)\big]<0. \]

{\bf Regime~II}: The unique optimal microstructure is a rank-I laminate.

This regime occurs when the function
\[ [0,\gamma^*]\ni\beta\mapsto\det\left[\left(d_2\alpha_1+d_1\alpha_2
-\beta T\right)^{-1}\left(\alpha_2(\ve_2^T-\ve)-
\alpha_1(\ve_1^T-\ve)\right)\right] \]
has a unique root (which we denote by $\beta_{II}$).

{\bf Regime~III}: There exist two optimal rank-II laminates.
This regime is present if the operator
$(d_2\alpha_1+d_1\alpha_2-\gamma^*T)$ is invertible and
\[ \det\left[\left(d_2\alpha_1+d_1\alpha_2-\gamma^*T\right)^{-1}
\left(\alpha_2(\ve_2^T-\ve)-\alpha_1(\ve_1^T-\ve)\right)\right]>0. \]

\begin{figure}[h!t]
\unitlength1cm
\begin{picture}(13.,3.2)
\if\Bilder y \put(+4.6,0.0){\psfig{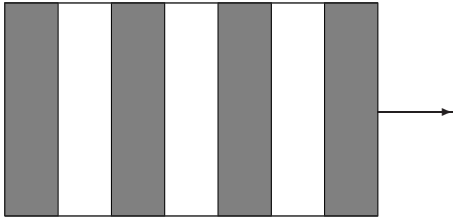}}\fi
\put(9.57,1.4){\vector(1,0){1}}
\end{picture}
\caption{\label{fig1}A two-phase rank-I laminate in two space dimensions with
corresponding normal vector. The strains are constant in the shaded and in the
unshaded regions. The volume fraction of both phases, $0.5$ in the picture,
is prescribed by the macroscopic parameter $d$.}
\end{figure}

\begin{figure}[ht]
\unitlength1cm
\begin{picture}(13.,3.2)
\if\Bilder y \put(+3.45,-6.1){\psfig{figure=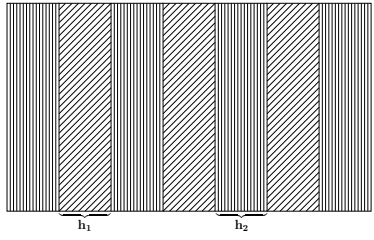,width=7.3cm}}\fi
\end{picture}
\caption{\label{fig2}A two-phase rank-II laminate in two space dimensions. The
widths $h_1$ and $h_2$ of the slabs should be much larger than the thickness
of the layers between the slab.}
\end{figure}

\vspace*{2mm}
For illustration, we visualize prototypes of rank-I and rank-II laminates
in Figures~\ref{fig1} and \ref{fig2}.
\vspace*{3mm}

To complete the definition (\ref{phidef}), set
\begin{align}
\alpha(\beta^*,d) &:= d_2\alpha_1+d_1\alpha_2-\beta^*T,\nn\\
e(\ve) &:= \alpha_2(\ve_2^T-\ve)-\alpha_1(\ve_1^T-\ve),\nn\\ 
\label{epsistardef}
\ve_i^*\equiv\ve_i^*(\beta^*,d,\ve) &:= \alpha^{-1}(\beta^*,d)
e_i(\beta^*,d,\ve),\\
e_1(\beta^*,d,\ve) &:= (\alpha_2-\beta^*T)\ve-d_2(\alpha_2\ve_2^T
-\alpha_1\ve_1^T),\nn\\
e_2(\beta^*,d,\ve) &:= (\alpha_1-\beta^*T)\ve+d_1(\alpha_2\ve_2^T
-\alpha_1\ve_1^T)\nn.
\end{align}
Hence $\ve_2^*-\ve_1^*=[\alpha(\beta^*,d)]^{-1}e(\ve)$.

Explicit computations of $\pad\whW$ and $\pae\whW$ are lengthy. We recall
the following results for the partial derivatives of $\whW$ in $D=2$ which are
proved in \cite{BC11}. Set
\begin{align*}
\sigma_i^* &:= \alpha_i(\ve_i^*-\ve_i^T),\\
\ovs^* &:= d_1\sigma_1^*+d_2\sigma_2^*.
\end{align*}

\begin{lemma}
\label{lem1}
Let $D=2$ and let $\alpha_i$ and $T$ commute. Then
\begin{multline}
\label{dWdeps}
\frac{\partial\whW}{\partial\ve}(d,\ve)=d_1\alpha_1\left(\ve_1^*-\ve_1^T\right)
+d_2\alpha_2\left(\ve_2^*-\ve_2^T\right)\\
+\begin{cases}
  \gamma^*d_1d_2\alpha^{-1}(\gamma^*,d)(\alpha_1-\alpha_2)T
  \left(\ve_2^*-\ve_1^*\right) & \text{ in Regime III} \\
  0 & \text{ else}.
\end{cases}
\end{multline}
Alternatively, in Regime~III,
\begin{multline}
\label{dWdepsIII}
\frac{\partial\whW}{\partial\ve}(d,\ve)
 = d_1(\alpha_2-\gamma^*T)\alpha^{-1}(\gamma^*,d)
   \alpha_1\left(\ve_1^*(\gamma^*,d,\ve)-\ve_1^T\right)\\
   +d_2(\alpha_1-\gamma^*T)\alpha^{-1}(\gamma^*,d)
   \alpha_2\left(\ve_2^*(\gamma^*,d,\ve)-\ve_2^T\right).
\end{multline}
Moreover it holds
\begin{equation}
\label{dWdchi}
\frac{\partial\whW}{\partial d}(d,\ve)
 =\ovs^*\!:\!\triangle\ve^*+W_1(\ve_1^*)-W_2(\ve_2^*)\\
+\begin{cases}
  0 & \text{ in Regimes 0 and I},\\
  \beta^*d_1d_2\frac{\partial\beta^*}{\partial d}
  \|\triangle\ve^*\|^2 & \text{ in Regime II},\\
  (d_1-d_2)\gamma^*\varphi(\triangle\ve^*) & \text{ in Regime III}.
\end{cases}
\end{equation}
\begin{align}
\frac{\partial\beta^*}{\partial d}
 &=\begin{cases}
\frac{\left(T(d_2\alpha_1+d_1\alpha_2-\beta T)^{-1}
(\alpha_2-\alpha_1)\triangle\ve^*\right):\triangle\ve^*}{\left(
(d_2\alpha_1+d_1\alpha_2-\beta T)^{-1}(T\triangle\ve^*)\right):
T\triangle\ve^*} & \text{in Regime II},\\
0 & \text{otherwise}.
\end{cases}
\label{betrus} \\[2mm]
\frac{\partial\beta^*}{\partial\ve} &=\begin{cases} \frac{1}{\left(
(d_2\alpha_1+d_1\alpha_2-\beta T)^{-1}(T\triangle\ve^*)\right):
T\triangle\ve^*}(\alpha_1-\alpha_2)\alpha^{-1}T\triangle\ve^* &
\text{in Regime II},\\
0 & \text{otherwise}.
\end{cases}\label{eq:beps}
\end{align}
\end{lemma}

\subsection{The scalar setting for $D=3$}
\label{secpoly}
The scalar setting in three dimensions is characterized by the ansatz
\begin{equation}
\label{uscal}
\bu(x_1,x_2,x_3):=\left(\!\!\begin{array}{c}x_1\\ x_2\\ \eta(x_1,x_2)\end{array}
\!\!\!\!\right)
\end{equation}
for the deformation, where $\eta$ is a scalar function (hence the name 'scalar'
theory, see \cite{BK97}, \cite{BS2}). Physically, (\ref{uscal}) corresponds to
anti-plane shear in the $x_3$-plane. Since
\[ \ve(\bu)=\left(\!\!\begin{array}{ccc} 1 & 0 & \partial_{x_1}\eta\\ 0 & 1 &
\partial_{x_2}\eta\\
\partial_{x_1}\eta &\partial_{x_2}\eta & 0 \end{array}\!\!\!\!\right), \]
$\nabla\eta$ determines the strain. This justifies to work with vectors
$f=\nabla\eta\in\R^2$, not with matrices $\ve(\bu)$. We replace (\ref{Wmin}) by
\begin{equation}
\label{Wmins}
\whW(d,f):=\!\inf_{\substack{\langle\td\rangle=d\\
\td\in\{0,1\}}}\;\inf_{\teta_{|\partial\Omega'}=f\cdot x}\slintOp\td_1
W_1(\nabla\teta)+\td_2W_2(\nabla\teta)\dx,\quad d\in[0,1],
\end{equation}
where the domain of integration $\Omega'$ is now two-dimensional. Therefore,
$\td\in BV(\Omega';\,\{0,1\})$, and it is well-understood from the context that
\[ \langle\td\rangle:=\frac{1}{|\Omega'|}\int_{\Omega'}\td(x)\dx \]
denotes the two-dimensional average. In addition, in (\ref{Wmins}) we adapted
the common notation for the scalar case and wrote $f$ instead of $\ve$ for the
strain. Finally, we set
\[ W_i(f):=\half\alpha_i(f-f_i^T)\cdot(f-f_i^T)+w_i \]
with $f_i^T$ the transformation strains in the scalar setting and $\alpha_i$
positive definite matrices.

\vspace*{2mm}
In \cite{CB2008}, partial results for $\whW$ in the non-scalar, i.e. vectorial
setting are available for $D=3$.
Yet, in general, the computation of $\pad\whW$, $\pae\whW$ for $D=3$
as required for assumption (\cA) remains currently open due to its complexity.
This is the reason why we restrict ourselves to a special three-dimensional
case. Our main result in this subsection is the following theorem that provides
us with an explicit formula for $\whW$ in the scalar case.
\begin{theorem}[Representation formula for $\whW$ and $D=3$ in the scalar
setting]
\label{theo5}
The functional $\whW$, given by (\ref{Wmins}), satisfies the explicit
representation formula
\begin{equation}
\label{whWsexp}
\whW(d,f)=d_1W_1(\nabla\eta_1^*)+d_2W_2(\nabla\eta_2^*)
\end{equation}
with
\begin{align}
\label{etasdef1}
\nabla\eta_1^*(d,f) &:= (d_2\alpha_1+d_1\alpha_2)^{-1}\big[\alpha_2f-d_2(
\alpha_2f_2^T-\alpha_1f_1^T)\big],\\
\label{etasdef2}
\nabla\eta_2^*(d,f) &:= (d_2\alpha_1+d_1\alpha_2)^{-1}\big[\alpha_1f+d_1(
\alpha_2f_2^T-\alpha_1f_1^T)\big].
\end{align}
\end{theorem}

\vspace*{2mm}
Theorem~\ref{theo5} states in particular that the scalar elastic theory in 3D is
directly related to the non-scalar geometrically linear elasticity theory in 2D
which was discussed in Subsection~\ref{subexpl2d}.

\vspace*{2mm}
{\bf Proof}: 
We apply the translation method, see for instance \cite{Che00}, \cite{Mil01}
for an overview. In general, the function $\phi$ which describes the translation
is only quasi-convex. Here, as a benefit of the scalar theory, we can pick
$\phi$ as a convex function, cf.\ \eqref{phipick}.

Starting from the inequality
\[ \phi(f)\le\inf_{\eta|\partial\Omega'=f\cdot x}\slintOp\phi(\nabla\eta)\dx, \]
which is satisfied by any quasi-convex function $\phi$, we obtain,
similar to the reasoning in \cite{CB2008} for the non-scalar case in $D=2$
dimensions,
\begin{align}
\whW(d,f)\ge\!\max_{\substack{\beta\ge0\\W_i-\beta\phi
\mbox{\tiny{ convex}}}}\hspace*{-2pt}
\min_{\substack{\td\in BV(\Omega';\{0,1\})\\\langle\td\rangle=d}}\;
\inf_{\eta|\partial\Omega'=f\cdot x}\!\Bigg[\!\! &
\slintOp\td_1(W_1-\beta\phi)(\nabla\eta)+d_2(W_2-\beta\phi)(\nabla\eta)\dx\nn\\
\label{star1}
& \hspace*{8pt} +\beta\phi(f)\!\Bigg].
\end{align}
Consider
\begin{equation}
\label{h1}
\nabla\eta_1(\td,\eta):=\frac{\slintOp\td_1\nabla\eta\dx}{\slintOp\td_1\dx},
\qquad
\nabla\eta_2(\td,\eta):=\frac{\slintOp\td_2\nabla\eta\dx}{\slintOp\td_2\dx}.
\end{equation}
Then
\begin{equation}
\label{eteq}
d_1\nabla\eta_1(\td,\eta)+d_2\nabla\eta_2(\td,\eta)=f.
\end{equation}
Next we apply Jensen's inequality which is possible since $W_i-\beta\phi$ is
convex for $i=1,2$. Furthermore we possibly enlarge the set of admissible
functions in the minimization over $\eta$ in \eqref{star1} by considering the
set of admissible functions $\nabla\eta_1, \nabla\eta_2\in\R^2$ with the
constraint \eqref{eteq}. Then \eqref{star1} becomes
\[ \whW(d,f)\ge\!\!\max_{\substack{\beta\ge0\\W_i-\beta\phi
\mbox{\tiny{ convex}}}}
\;\,\min_{\substack{\nabla\eta_1,\nabla\eta_2\in\R^2\\ d_1\nabla\eta_1
+d_2\nabla\eta_2=f}}
\Big\{d_1W_1(\nabla\eta_1)+d_2W_2(\nabla\eta_2)+\beta\phi(f)\Big\}. \]

If $D=2$, $\phi$ is chosen to be $-\det(\ve(\bu))$. Here we set analogously
\begin{equation}
\label{phipick}
\phi(\nabla\eta):=\big(\partial_x\eta\big)^2
+\big(\partial_y\eta\big)^2\ge0.
\end{equation}
Since $\phi$ is quadratic, there exists a unique linear operator
$T:\R^2\to\R^2$ such that $\phi(f)=\half Tf\!:\!f$, here simply $Tf=2f$.
It holds
\begin{align*}
& \hspace*{-20pt} \sum_{i=1}^2d_i(W_i-\beta\phi)(\nabla\eta_i)+\beta\phi(f)\\
& =\sum_{i=1}^2
d_iW_i(\nabla\eta_i)+\beta\big[\underbrace{\phi(d_1\nabla\eta_1+d_2\nabla\eta_2)
-d_1\phi(\nabla\eta_1)-d_2\phi(\nabla\eta_2)}_{=:R}\big].
\end{align*}
By a direct computation the remainder term $R$ can be rewritten,
\begin{eqnarray*}
2R &=& T(d_1\nabla\eta_1+d_2\nabla\eta_2)\!:\!(d_1\nabla\eta_1+d_2\nabla\eta_2)
-d_1T\nabla\eta_1\!:\!\nabla\eta_1-d_2T\nabla\eta_2\!:\!\nabla\eta_2\\
&=& d_1d_2\big[-T\nabla\eta_1\!:\!\nabla\eta_1-T\nabla\eta_2\!:\!\nabla\eta_2
+T\nabla\eta_1\!:\!\nabla\eta_2+T\nabla\eta_2\!:\!\nabla\eta_1\big]\\
&=& -d_1d_2T(\nabla\eta_2-\nabla\eta_1)\!:\!(\nabla\eta_2-\nabla\eta_1)\\
&=& -2d_1d_2\phi(\nabla\eta_2-\nabla\eta_1).
\end{eqnarray*}
So we have found
\begin{equation}
\label{h2}
\whW(d,f)\ge\!\!\max_{\substack{\beta\ge0\\W_i-\beta\phi
\mbox{\tiny{ convex}}}}
\;\,\min_{\substack{\nabla\eta_1,\nabla\eta_2\in\R^2\\ d_1\nabla\eta_1
+d_2\nabla\eta_2=f}}\!\!
\Big\{d_1W_1(\nabla\eta_1)+d_2W_2(\nabla\eta_2)-\beta d_1d_2
\phi(\nabla\eta_2\!-\!\nabla\eta_1)\Big\}.
\end{equation}
Next we compute the optimal strains $\nabla\eta_1^*$, $\nabla\eta_2^*$.
After differentiating the argument on the right in (\ref{h2}) for fixed
$\beta$, we obtain
\begin{equation}
\label{nosjump}
\alpha_1(\nabla\eta_1^*-f_1^T)-\alpha_2(\nabla\eta_2^*-f_2^T)+\beta T
(\nabla\eta_2^*-\nabla\eta_1^*)=0.
\end{equation}
Using the constraint $d_1\nabla\eta_1^*+d_2\nabla\eta_2^*=f$, after
rearrangement this yields the formulas
\begin{align*}
(d_1\alpha_2+d_2\alpha_1-\beta T)\nabla\eta_1^* &= (\alpha_2-\beta T)f
-d_2(\alpha_2f_2^T-\alpha_1f_1^T),\\
(d_1\alpha_2+d_2\alpha_1-\beta T)\nabla\eta_2^* &= (\alpha_1-\beta T)f
+d_1(\alpha_2f_2^T-\alpha_1f_1^T).
\end{align*}
Setting $\beta=0$ (see below), this coincides with (\ref{etasdef1}),
(\ref{etasdef2}).

The maximum over $\beta$ in (\ref{h2}) is attained at $\beta^*=0$ as
$d_1d_2\phi(\nabla\eta_2^*-\nabla\eta_1^*)\ge0$. Hence, the translation
is trivial in this setting. We obtain
\begin{equation}
\label{sest2}
\whW(d,f)\ge d_1W_1(\nabla\eta_1^*)+d_2W_2(\nabla\eta_2^*)=:\whW_-(d,f).
\end{equation}
It remains to estimate $\whW$ from above which amounts to showing that
$\nabla\eta_1^*$, $\nabla\eta_2^*$ yield an optimal microstructure.
Plugging in any microstructure $\td$ with $\langle\td\rangle=d$, from the
definition (\ref{Wmins}) of $\whW$, we get the upper bound
\[ \whW(d,f)\le\inf_{\teta_{|\partial\Omega'}=f\cdot x}\slintOp\td_1
W_1(\nabla\eta)+\td_2W_2(\nabla\eta)\dx=:\whW_+(d,f). \]
As the domain of integration $\Omega'$ in the scalar case with $D=3$ is
two-dimensional, finding the optimal microstructure subsumes to the non-scalar
case in $D=2$ dimensions. Hence, as can be seen from (\ref{bsdef}), $\beta^*=0$
occurs either for Regime~0 (no microstructure) or Regime~I. The latter is
equivalent to $\phi(\ve_2^*-\ve_1^*)>0$. From a well-known argument, see
\cite[Lemma~4.1]{Kohn91},
this implies that $\ve_1^*$, $\ve_2^*$ are compatible.
Besides this compatibility condition, any optimal microstructure must also
satisfy the equilibrium condition
\begin{equation}
\label{state}
[[\sigma]]\vec{n}=0
\end{equation}
which states that the jump of the stress in the normal direction of  the
laminate must vanish. Equation~(\ref{state}) originates from the Euler-Lagrange
equation of the variational problem. Here, it is automatically satisfied since,
due to (\ref{nosjump}), $[[\sigma]]\vec{n}=\sigma_2^*-\sigma_1^*=0$. 

So there exists a unique rank-I lamination microstructure
(unique up to a sign $\pm1$ of direction of the laminate) that connects
$\ve_1^*$ and $\ve_2^*$, and the strain of phase $i$ is $\nabla\eta_i^*$,
$i=1,2$. This microstructure is optimal, $\whW_+(d,f)=\whW_-(d,f)$,
which proves (\ref{whWsexp}). \qed

\subsection{The case $D=1$}
\label{secd1}
The methods of the previous subsection also apply to derive rigorously an
explicit expression for $\whW$ in $D=1$. Here, quasi-convex functions are
always convex, and the vector and scalar settings coincide. However, we use the
notation in analogy to the vector case.

\begin{theorem}[Representation formula for $\whW$ in $D=1$]
\label{theoa}
For $D=1$, the functional $\whW$, given by (\ref{Wmin}), satisfies the explicit
representation formula
\begin{equation}
\label{Wh1d}
\whW(d,\ve)=dW_1(\ve_1^*)+(1-d)W_2(\ve_2^*)
\end{equation}
with
\begin{align}
\label{e1db}
\ve_1^*(d,\ve) &= \frac{\alpha_2(\ve-d_2\ve_2^T)+d_2\alpha_1\ve_1^T}
{d_2\alpha_1+d_1\alpha_2},\\
\label{e1dc}
\ve_2^*(d_1,\ve) &= \frac{\alpha_1(\ve-d_1\ve_1^T)+d_1\alpha_2\ve_2^T}
{d_2\alpha_1+d_1\alpha_2}.
\end{align}
The partial derivatives of $\whW$ are given by
\begin{align}
\label{dWdeps1d}
\frac{\partial\whW}{\partial\ve}(d,\ve) &= \frac{\alpha_1\alpha_2\big[
d_1(\ve-\ve_1^T)+d_2(\ve-\ve_2^T)\big]}{d_2\alpha_1+d_1\alpha_2},\\
\frac{\partial\whW}{\partial d}(d,\ve) &= W_1(\ve_1^*)-W_2(\ve_2^*)
+\frac{\alpha_1\alpha_2}{(d_2\alpha_1+d_1\alpha_2)^2}\Big[
(\alpha_1-\alpha_2)\ve^2+d_1\alpha_1(\ve_1^T)^2-d_2\alpha_2(\ve_2^T)^2\nn\\
& \hspace*{170pt}
+\big(\alpha_2\ve_2^T\!-\!\alpha_1\ve_1^T
\!+\!(\alpha_2\!-\!\alpha_1)(d_1\ve_1^T\!+\!d_2\ve_2^T)\big)\ve\nn\\
& \hspace*{170pt}
\label{dWdd1d} +(d_2\alpha_1-d_1\alpha_2)\ve_1^T\ve_2^T\Big].
\end{align}
\end{theorem}

\vspace*{2mm}
{\bf Proof}:
We apply the same methods as in the proof of Theorem~\ref{theo5}.

First, we estimate $W_1$, $W_2$ from below by their convexification. Similar to
(\ref{h2}), we have for any convex function $\varphi$ and $\beta\ge0$,
\[ \whW(d,\ve)\ge\max_{\substack{\beta\ge0\\ W_i-\beta\varphi\;
\mbox{{\tiny convex}}}}\;\min_{\substack{\ve_1,\ve_2\in\R^D\\
d\ve_1+(1-d)\ve_2=\ve}}\big\{d(W_1-\beta\varphi)(\ve_1)
+(1-d)(W_2-\beta\varphi)(\ve_2)+\beta\varphi(\ve)\big\}. \]

Secondly, we choose the optimal $\beta$ and $\varphi$. The later is an art when
$D>1$. With arguments identical to the scalar case in $D=3$, it holds
$\beta\equiv0$, leaving us with the estimate
\[ \whW(d,\ve)\ge\min_{\substack{\ve_1,\ve_2\in\R\\ d\ve_1+(1-d)\ve_2=\ve}}
\big\{dW_1(\ve_1)+(1-d)W_2(\ve_2)\big\}=:W_-(d,\ve). \]

Next we compute the optimal strains $\ve_1$, $\ve_2$ on the right. For $d=1$,
only $\ve_1=\ve$ is admissible and hence optimal. Now, let $0\le d<1$.
We can then resolve the constraint by setting
$\ve_2=\frac{\ve-d\ve_1}{1-d}$ and need to calculate
\begin{equation}
\label{star2}
\min_{\ve_1\in\R}\Bigg\{dW_1(\ve_1)+(1-d)W_2\Big(\frac{\ve-d\ve_1}{1-d}\Big)
\Bigg\}.
\end{equation}
After differentiation with respect to $\ve_1$ and using (\ref{Widef})
for the optimal $\ve_1^*$ in (\ref{star2}), we are left with
\[ \alpha_1(\ve_1^*-\ve_1^T)=\alpha_2\Big(
\frac{\ve-d\ve_1^*}{1-d}-\ve_2^T\Big). \] 
Rearrangement and resolution of $d\ve_1^*+(1-d)\ve_2^*=\ve$ gives the
formulas (\ref{e1db}), (\ref{e1dc}) of the optimal strains. Since
$\ve_1^*(1,\ve)=\ve$, these formulas also state the correct solution when $d=1$.

In the final step, due to the definition of $\whW$, it holds
\[ \whW(d,\ve)\le\inf_{\bu|\partial\Omega=\ve x}\slintO\td W_1(\ve(\bu))+(1-\td)
W_2(\ve(\bu))\dx=:W_+(d,\ve) \]
and $\td\in BV(\Omega;\,\{0,1\})$ may represent any microstructure with
$\langle\td\rangle=d$. The explicit
construction of the optimal microstructure for which the upper bound and the
lower bound coincide,
\begin{equation}
\label{con1} W_-(d,\ve)=\whW(d,\ve)=W_+(d,\ve),
\end{equation}
is again as in the proof of Theorem~\ref{theo5}. 
This leads to (\ref{Wh1d}).

Finally, the verification of (\ref{dWdeps1d}), (\ref{dWdd1d}) follows from
\begin{equation}
\label{depss1D1}
\frac{\partial\ve_1^*}{\partial\ve}(d,\ve)=\frac{\alpha_2}
{d_2\alpha_1+d_1\alpha_2},\qquad
\frac{\partial\ve_2^*}{\partial\ve}(d,\ve)=\frac{\alpha_1}
{d_2\alpha_1+d_1\alpha_2}
\end{equation}
and
\begin{align}
\label{depss1D2a}
\frac{\partial\ve_1^*}{\partial d} &= \frac{\alpha_2\big[\alpha_2\ve_2^T
-\alpha_1\ve_1^T+(\alpha_1-\alpha_2)\ve\big]}{(d_2\alpha_1+d_1\alpha_2)^2},\\
\label{depss1D2b}
\frac{\partial\ve_2^*}{\partial d} &= \frac{\alpha_1\big[\alpha_2\ve_2^T
-\alpha_1\ve_1^T+(\alpha_1-\alpha_2)\ve\big]}{(d_2\alpha_1+d_1\alpha_2)^2}
\end{align}
which all can be verified by elementary computations.

Using the relationship $d_1\ve_1^*+d_2\ve_2^*=\ve$ we find (\ref{dWdeps1d}) and
finally
\begin{align}
\label{star3}
\frac{\partial\whW}{\partial d}(d,\ve) &= W_1(\ve_1^*)-W_2(\ve_2^*)
+d_1W_1'(\ve_1^*)\frac{\partial\ve_1^*}{\partial d}
+d_2W_2'(\ve_2^*)\frac{\partial\ve_2^*}{\partial d},\\
&= W_1(\ve_1^*)-W_2(\ve_2^*)\nn\\
& \quad +\frac{\alpha_1\alpha_2\big(\alpha_2\ve_2^T-\alpha_1\ve_1^T+(\alpha_1
-\alpha_2)\ve\big)\big[d_1(\ve_1^*-\ve_1^T)+d_2(\ve_2^*-\ve_2^T)\big]}
{(d_2\alpha_1+d_1\alpha_2)^2}.\nn
\end{align}
Using
\begin{equation}
\label{star4}
d_1(\ve_1^*-\ve_1^T)+d_2(\ve_2^*-\ve_2^T)=\ve-d_2\ve_2^T-d_1\ve_1^T,
\end{equation}
this simplifies to identity (\ref{dWdd1d}). \qed

\section{Existence results for the AC-CH system with microstructural energy
densities}
\label{secexam}
Above we have already remarked that the energy $W_{\mathrm{lin}}$ satisfies
assumption (\cA) and thus Theorem~\ref{theo1} applies. In this section
we investigate whether the results transfer to the geometrically linear
theory of elasticity stated in the previous section.

The following statement is a consequence of Theorem~\ref{theo1}.

\begin{theorem}[Existence of solutions for geometrically linear elastic energy]
\label{thmcor1}
Let $\cW=\whW$ be a function on $[0,1]\times\R^{D\times D}_\mathrm{sym}$.
Assume that the mobility tensor $M$ be positive definite \cb and continuous
for all $a$, $b$ satisfying (\ref{abrange}). \cn Moreover, let
$\psi$ be given by (\ref{psidef}) and the initial data $(a_0,b_0)$ satisfy
(\ref{abrange}) \cb and (\ref{W0bound}). \cn

In $D\!=\!3$, let $\whW$ be given by (\ref{whWsexp}) corresponding to the
scalar setting. If $D=2$, let
\begin{enumerate}
 \item[(\cC1)] $\beta^*(d,\ve)$ be independent of $\ve$, and
 \item[(\cC2)] $\alpha_i$ and $T$ commute whenever $\beta^*\in
\{\gamma^*, \beta_{II}\}$.
\end{enumerate}
If $D=1$, let $\whW$ be given by (\ref{Wh1d}).
Then there exists a solution $(a,b,\bu)$ to (\ref{ACH1})--(\ref{ACH3})
that satisfies

(i) $a,b\in C^{0,\frac14}\left([0,{\cal T}];\,L^2(\Omega)\right)$,\\
(ii) \cb$\pat a\in L^2(0,{\cal T};\,H^1(\Omega)^*)$, \cn
$\pat b\in L^2(\OT)$,\\
(iii) $\bu\in L^\infty\left(0,{\cal T};\,H^1(\Omega;\,\R^D)\right)$,\\
(iv) The feasible parameter range of $(a,b)$ is given by (\ref{abrange}).
\end{theorem}

\vspace*{2mm}
{\bf Proof}: We first consider the case $D<3$.
The idea is to apply Theorem~\ref{theo1}. Therefore we 
\begin{enumerate}
\item verify that $\whW$ satisfies assumption (\cA) for $d\in[0,1]$,
\item extend $\whW$ to $\R\times\R^{D\times D}_\mathrm{sym}$ in such a way that
(\cA) is fulfilled for all $d\in\R$.
\end{enumerate}
Then Theorem~\ref{theo1} can be applied and we obtain by part (iv)
therein that effectively $d\in [0,1]$ in the system
(\ref{ACH1})--(\ref{ACH3}). Hence Theorem~\ref{thmcor1} follows with
$\cW=\whW$ as asserted.

We start with the second step and give the proofs of the former in dimensions
$D=1$, $D=2$, and the scalar case in $D=3$ thereafter. We introduce the
following extension $\whW_0(\cdot,\ve):\R\to\R$ of $\whW(\cdot,\ve)$ by
\begin{equation}
\label{Wext}
\whW_0(d,\ve):=\left\{\begin{array}{ll}
\whW(d,\ve), & 0\le d\le 1,\\
\rho_1(d,\ve), & 0>d>-1,\\
-d+1+\whW(0,\ve), & d\le-1,\\
\rho_2(d,\ve), & 1<d<2,\\
d-1+\whW(1,\ve), &d\ge2 \end{array}\right.
\end{equation}
for functions $\rho_1$, $\rho_2$ that we may choose such that
$\whW_0(\cdot,\cdot)\in C^1(\R\times\R^{D\times D}_\mathrm{sym})$ and that
$\whW_0$ fulfils (\cA1)--(\cA2) for $d<0$ and $d>1$. 

We further mention that the function $\whW_0$ defined in (\ref{Wext}) is
positive and coercive for $|d|\to\infty$. We note that this construction is
required for the first part of the proof of Theorem~\ref{theo1}, where
polynomial approximations of $\psi$ are constructed, where $d\notin[0,1]$ may
occur.

It remains to show that $\whW$ fulfils (\cA) for $d\in [0,1]$, which we
prove independently for $D=1$, $D=2$, and the scalar case in $D=3$.

(i) {\it Verification of (\cA) for $D=1$}:

From~\eqref{Widef}, (\ref{Wh1d}), \eqref{e1db} and \eqref{e1dc}, it follows
$\whW\in C^1(\R\times\R)$.

From (\ref{dWdeps1d}), we obtain
\[ \left(\pae\whW(d,\ve_2)-\pae\whW(d,\ve_1)\right)(\ve_2-\ve_1)
=\frac{\alpha_1\alpha_2}{d_2\alpha_1+d_1\alpha_2}\sum_{i=1}^2
d_i(\ve_i^*(d,\ve_2)-\ve_i^*(d,\ve_1))(\ve_2-\ve_1). \]
The condition (\cA1), restricted to $d\in[0,1]$, follows, since by
\eqref{e1db}, \eqref{e1dc} for $i=1,2$,
\[ \left(\ve_i^*(d,\ve_2)-\ve_i^*(d,\ve_1)\right)(\ve_2-\ve_1)\ge
\tfrac{\min\{\alpha_1,\alpha_2\} }{\max\{\alpha_1,\alpha_2\}}|\ve_2-\ve_1|^2, \]
where $\tfrac{\min\{\alpha_1,\alpha_2\}}{\max\{\alpha_1,\alpha_2\}}>0$ by
assumption. From \eqref{e1db}, \eqref{e1dc}, for $i=1,2$ and $d\in[0,1]$ we have
by similar arguments
\begin{equation}
\label{visest}
|\ve_i^*(d,\ve)|\le c\big(|d|+|\ve|+1\big).
\end{equation}
This leads to (\cA2)$_1$, since by (\ref{Wh1d}), \eqref{visest},
\begin{align*}
|\whW(d,\ve)| &\le |W_1(\ve_1^*)|+|W_2(\ve_2^*)|
\le c\left(|\ve_1^*|^2+|\ve_2^*|^2+1\right)\\
&\le c\left(|d|^2+|\ve|^2+1\right).
\end{align*}

For $d\in[0,1]$ we easily compute, for $i=1,2$,
\[ |\pad\ve_i^*(d,\ve)|\le c\left(|\ve|+1\right). \]

We use this estimate in (\ref{star3}), to obtain by \eqref{visest}
\begin{align*}
|\pad\whW(d,\ve)| &\le c\left(|\ve_1^*(d,\ve)|^2+|\ve_2^*(d,\ve)|^2 +1
+(|\ve_1^*(d,\ve)|+|\ve_2^*(d,\ve)|)(|\ve|+1)\right)\\
&\le c\left(|d|^2+|\ve|^2+1\right)
\end{align*}
which is (\cA2)$_2$ restricted to $d\in[0,1]$.

For $d\in[0,1]$ we derive with the help of (\ref{dWdeps1d}), 
(\ref{star4}), \eqref{visest}
\begin{align*}
|\pae\whW(d,\ve)| &\le\left|
\tfrac{\alpha_1\alpha_2}{\min\{\alpha_1,\alpha_2\}}\right|
\left(|\ve-\ve_1^T|+|\ve-\ve_2^T|\right)\\
&\le \max\{\alpha_1,\alpha_2\}\; c\left(|\ve|+1\right)\\ 
&\le c\left(|\ve|+1\right),
\end{align*}
showing the validity of (\cA2)$_3$ restricted to $d\in[0,1]$.

\vspace*{3mm}
(ii) {\it Verification of (\cA) for $D=2$}:

The regularity of $\whW$ required for (\cA) follows directly from \eqref{Wdef}
and the related definitions. 

The constant $\gamma^*$ solely depends on $\alpha_1$, $\alpha_2$ since
$\gamma^*=\min\{\gamma_1,\gamma_2\}$ and $\gamma_i$ is the reciprocal of the
largest eigenvalue of $\alpha_i^{-1/2}T\alpha_i^{-1/2}$, see \eqref{gsdef}.
Since $\beta^*\in[0,\gamma^*]$,
\begin{equation}
\label{bsest}
|\beta^*(d,\ve)|\le c
\end{equation}
for a constant $c$ independent of $d$ and $\ve$.

From (\ref{bsest}), (\ref{epsistardef}) we obtain for $d\in[0,1]$
\begin{equation}
\label{e1D2} 
|\ve_i^*(\beta^*,d,\ve)|\le c\big(|d|+|\ve|+1).
\end{equation}
With~\eqref{Widef} and \eqref{Wdef} this shows for $d\in[0,1]$
\begin{align*}
\left|\beta^*d_1d_2\det(\ve_2^*(d,\ve)-\ve_1^*(d,\ve))\right|
 &\le c\left(|\ve_1^*(d,\ve)|^2+|\ve_2^*(d,\ve)|^2+1\right)\\
 &\le c\left(|d|^2+|\ve|^2+1\right).
\end{align*}
From this we easily verify (\cA2)$_1$ restricted to $d\in[0,1]$, as the
first two terms in (\ref{Wdef}) can be estimated as in the one-dimensional case.

In order to show (\cA2)$_2$, we apply \eqref{dWdchi} and (\ref{visest}) to find
\begin{align*}
|\ovs^*(d,\ve)| &= \left|\sum_{i=1}^2d_i\alpha_i(\ve_i^*-\ve_i^T)\right|\\
&\le c\left(|d|+|\ve|+1\right),\\
|\triangle\ve^*(d,\ve)| &= |\ve_2^*-\ve_1^*|\\
 &\le c\left(|d|+|\ve|+1\right),\\
|\varphi(\triangle\ve^*(d,\ve))| &= |\det(\triangle\ve^*(d,\ve))|
\le|\ve_2^*(d,\ve)-\ve_1^*(d,\ve)|^2\\
&\le c\left(|d|^2+|\ve|^2+1\right).
\end{align*}
Finally, from (\ref{betrus}) and the fact that $T$ is an isometry w.r.t. the
Frobenius norm,
\[ \left|\frac{\partial\beta^*}{\partial d}(d,\ve)\right|\le c. \]
The terms $W_i(\ve_i^*(d,\ve))$, $i=1,2$ remaining in \eqref{dWdchi} can be
estimated as in the one-dimensional case. This verifies (\cA2)$_2$
restricted to $d\in[0,1]$.

Now we prove (\cA2)$_3$. For the regimes 0, I and II this follows directly from
(\ref{visest}). When in Regime~III, we start
from~\eqref{dWdepsIII} and find with~\eqref{e1D2} for $d\in[0,1]$
\[ |d_1(\alpha_2-\gamma^*T)\alpha^{-1}\alpha_1(\ve_1^*-\ve_1^T)|
\le c\left(|\ve_1^*|+1\right)\le c\left(|d|+|\ve|+1\right). \]

The validation of (\cA1) relies on \eqref{dWdeps}.
By assumption (i) of the theorem, $\frac{\partial\beta^*}{\partial\ve}=0$
and, as $d\in[0,1]$ is fixed in (\cA1), $\alpha^{-1}$ is a constant tensor.
Assumption~(\cA1) then follows from
\[ \left(\ve_1^*(d,\ve_2)-\ve_1^*(d,\ve_1)\right):(\ve_2-\ve_1)
 =\alpha^{-1}(\alpha_2-\beta^*T)(\ve_2-\ve_1):(\ve_2-\ve_1), \]
a similar equality for
$\big(\ve_2^*(d,\ve_2)-\ve_2^*(d,\ve_1)\big)\!:\!(\ve_2-\ve_1)$, and the
positive definiteness of $\alpha^{-1}$ and $(\alpha_i-\beta^*T)$,
$i=1,2$. Indeed, the latter is equivalent to
$\Id-\beta^*\alpha_i^{-1/2}T\alpha_i^{-1/2}$ positive definite, and
since $\beta^*\le\gamma^*$, this is implied by all eigenvalues of
$\Id-\gamma^*\alpha_i^{-1/2}T\alpha_i^{-1/2}$ being positive. However, from
(\ref{gsdef}) and the definition of $\gamma_i$, this is true.

\vspace*{2mm}
For the scalar case in $D=3$, comparing (\ref{Wdef}) and (\ref{whWsexp})
and since the partial derivatives coincide, the theory subsumes to the
geometrically linear
theory for $D=2$ and $\beta^*=0$. Hence we can proceed as above and show that
$\whW$ given by (\ref{whWsexp}) satisfies (\cA). The assumptions (\cC1), (\cC2)
are not needed here because $\beta^*=0$. \qed

\vspace*{2mm}
The following theorem states the uniqueness of weak solutions for the system
with geometrically linear theory.
\begin{theorem}[Uniqueness of weak solutions]
\label{theo4}
Let $\cW=\whW$ be given by (\ref{whWsexp}) for $D=3$ (scalar setting), by
(\ref{Wdef}) for $D=2$, or by (\ref{Wh1d}) for $D=1$.
Assume further that $M\equiv1$ and the elastic moduli of the two phases are
equal, i.e. $\alpha_1=\alpha_2$. For $D=2$, let $\det(\ve_2^T-\ve_1^T)\le0$.
Then the solution $(a,b,\bu)$ in
Theorem~\ref{thmcor1} is unique in the spaces stated there.
\end{theorem}

\vspace*{2mm}
{\bf Proof}: (i) {\it The case $D=1$.}

For $\alpha_1=\alpha_2$, the equations (\ref{e1db}), (\ref{e1dc}) read
\[ \ve_1^*(d,\ve)=\ve+d_2(\ve_1^T-\ve_2^T),\qquad
\ve_2^*(d,\ve)=\ve+d_1(\ve_2^T-\ve_1^T) \]
which implies $\ve_2^*-\ve_2^T=\ve_1^*-\ve_1^T$ leading to
$W_1(\ve_1^*)-w_1=W_2(\ve_2^*)-w_2$. This yields
\[ \whW(d,\ve)=\frac{\alpha_1}{2}\big(\ve_1^*-\ve_1^T\big)^2+d_1w_1+d_2w_2
=\frac{\alpha_1}{2}\big(\ve+d_1(\ve_2^T-\ve_1^T)-\ve_2^T\big)^2+d_1w_1+d_2w_2.
\]
In this special case, the partial derivatives of $\whW$ are
\begin{align*}
\frac{\partial\whW}{\partial d}(d,\ve) &= \alpha_1\big(\ve+d_1(\ve_2^T-\ve_1^T)
-\ve_2^T\big)(\ve_2^T-\ve_1^T)+w_1-w_2,\\
\frac{\partial\whW}{\partial\ve}(d,\ve) &= \alpha_1\big(\ve+d_1(\ve_2^T-\ve_1^T)
-\ve_2^T\big).
\end{align*}
Now we pass through the steps in the proof of Theorem~\ref{theo2} and make
the necessary modifications.
Let again $(a^k,b^k,u^k)$, $k=1,2$ be two solutions to
(\ref{ACH1})--(\ref{ACH3}) and set $a:=a^1-a^2$, $b:=b^1-b^2$, $d^k:=a^k+b^k$,
$\ve^k:=\ve(u^k)$ for $k=1,2$ and $d:=d^1-d^2=a+b$, $\ve:=\ve^1-\ve^2$.

The equations (\ref{test2}), (\ref{test3}) now read
\begin{align}
& \ioT\left[\pat b\eta+\lambda\nabla b\cdot\nabla\eta+\alpha_1\big(\ve+d(\ve_2^T
-\ve_1^T)\big)(\ve_2^T-\ve_1^T)\eta\right]\dx\dt\nn\\
\label{2test2}
& \quad =\!\ioT\left[\frac{\theta}{2}\left(g'(a^2\!+\!b^2)-g'(a^1\!+\!b^1)
+g'(a^1\!-\!b^1)-g'(a^2\!-\!b^2)\right)\eta+2\kappa_2b\eta\right]\dx\dt,\\
\label{2test3}
& \ioto\alpha_1\big(\ve+d(\ve_2^T-\ve_1^T)\big)\ve(u)\dx\dt=0.
\end{align}
Analogous to (\ref{step2}) we find
\begin{align}
& \hspace*{-20pt} \frac{\lambda}{2}\|a(t_0)\|_M^2+\ioto\lambda|\nabla a|^2
+\alpha_1\big(\ve+d(\ve_2^T-\ve_1^T)\big)a(\ve_2^T-\ve_1^T))\dx\dt\le\nn\\
\label{2step4}
& \ioto2\kappa_1a^2+\frac{\theta}{2}\Big[\big|g'(a^1\!+\!b^1)
-g'(a^2\!+\!b^2)\big|+\big|g'(a^1\!-\!b^1)-g'(a^2\!-\!b^2)\big|\Big]|a|\dx\dt.
\end{align}
Choosing $\eta:=(b^1-b^2){\cal X}_{(0,t_0)}$ in (\ref{2test2}) and adding
the resulting equation to (\ref{2step4}), we obtain
\begin{align}
\label{wnote}
& \frac{\lambda}{2}\|a(t_0)\|_M^2+\half\|b(t_0)\|_{L^2}
+\ioto\Big[\lambda\big(|\nabla a|^2+|\nabla b|^2\big)+\alpha_1\big(\ve+
d(\ve_2^T-\ve_1^T)\big)^2\Big]\dx\dt\\
& \le\!\ioto\!2(\kappa_1|a|^2\!+\!\kappa_2|b|^2)
+\!\!\ioto\!\!\frac{\theta}{2}\Big[\big|g'(a^1\!\!+\!\!b^1)
\!-\!g'(a^2\!\!+\!\!b^2)\big|\!+\!\big|g'(a^1\!\!-\!\!b^1)\!-\!g'(a^2\!\!-\!\!
b^2)\big|\Big](|a|\!+\!|b|)\dx\dt\nn.
\end{align}
It is noteworthy that we are only able to estimate 
$\int_{\Omega_{\cal T}}\alpha_1\big(\ve+d(\ve_2^T-\ve_1^T)\big)^2\dx\dt$ which
differs from the total mechanical energy
$\int_{\Omega_{\cal T}}\whW(d,\ve)\dx\dt$.

With the Lipschitz continuity of $g'$, and by applying the inequalities of
Young and Gronwall to (\ref{wnote}), we obtain as in the proof of
Theorem~\ref{theo2}
\[ \ioT\alpha_1\big(\ve(u)\big)^2\dx\dt=0. \]
With Korn's inequality, $u\equiv0$ in $\OT$.

(ii) {\it The scalar case in $D=3$.}

By Theorem~\ref{theo5}, for $\alpha_1=\alpha_2$,
\begin{align*}
\nabla\eta_1^*(d,f) &= f-d_2(f_2^T-f_1^T),\\
\nabla\eta_2^*(d,f) &= f+d_1(f_2^T-f_1^T).
\end{align*}
Consequently, $\nabla\eta_1^*-f_1^T=\nabla\eta_2^*-f_2^T$ leading to
$W_1(\nabla\eta_1^*)-w_1=W_2(\nabla\eta_2^*)-w_2$. Since $\beta^*=0$, the
formula for $\whW$ has now the same structure as for $D=1$ and the proof
follows exactly as in the first part (i).

\vspace*{2mm}
(iii) {\it The case $D=2$.}

By assumption, $\phi(\beta,d,\ve)>0$, which results in
$\beta^*\equiv0$ in $\OT$.  In addition, for $\alpha_1=\alpha_2$, by
Eqn.~(\ref{epsistardef}),
\begin{align*}
\ve_1^*(d,\ve) &= \alpha_1^{-1}\big((\alpha_1-\beta^*T)\ve
-d_2\alpha_1(\ve_2^T-\ve_1^T)\big),\\
\ve_2^*(d,\ve) &= \alpha_1^{-1}\big((\alpha_1-\beta^*T)\ve
+d_1\alpha_1(\ve_2^T-\ve_1^T)\big).
\end{align*}
Consequently,
\[ \phi(\beta,d,\ve)=-\det((\ve_2^*-\ve_1^*)(d,\ve))=-\det(\ve_2^T-\ve_1^T). \]

So again $\ve_1^*-\ve_1^T=\ve_2^*-\ve_2^T$, $W_1(\ve_1^*)-w_1=W_2(\ve_2^*)-w_2$,
and the proof can be carried out as before. \qed

\section{Conclusion}
\label{secconclude}
In this article we derived extensions of the Allen-Cahn/Cahn-Hilliard system
to elastic materials. In particular we included
(i) the linear elastic energy derived by Eshelby, (ii) a geometrically linear
theory of elasticity for $D\le2$, and (iii) a geometrically linear theory of
elasticity for $D=3$ in the scalar setting. For all three cases we showed in a
mathematically rigorous way the existence and the uniqueness of weak solutions,
asserting the correctness of our approach. As a future goal, it is desirable to
generalize the proof to the non-scalar, i.e. vectorial case in $D=3$ and to
$D\ge3$.

The generalized AC-CH models contain as special cases both the Allen-Cahn
\cite{BW05} and the Cahn-Hilliard equation \cite{Harald1} with linear
elasticity. In our work the existence and uniqueness of weak solutions to the
AC-CH model with linear elasticity follows in a straightforward way since
Assumption~(\cA) required in the existence result Theorem~\ref{theo1} can easily
be checked to hold true for the linear energy. This had not been shown before.

We point out that for the new cases (ii) and (iii), the formulas collected in
Section~\ref{secexpl} for the specific energy $\whW$ are essential for any
numerical studies of the AC-CH models extended to microstructure. For related
investigations and numerical methods we refer the interested reader to
\cite{BC11} and \cite{B08}.

Besides the limiting assumption of constant temperature (cf.\ Introduction),
the most important pending restriction is the postulation of small strain,
included in (\ref{vedef}). That is, it would be desirable to combine the AC-CH
system with energies within the framework of the geometrically non-linear
theory of elasticity. The technical problems of a large strain theory are
striking, see, e.g., \cite{Kaushik,BK97,CB2008,KN00},
and it is in many cases not known how to compute explicit formulas
for the relaxed energy functionals.

Similar mathematical difficulties arise when one wants to extend the AC-CH
systems to elastic materials which are of polycrystalline structure. This is
an open problem for future research. For related modelling aspects we refer to
\cite{BS}, for related mathematical aspects to \cite{BK97,BS2}.

\section*{Acknowledgment}
We have worked on this project while we were employed or were
visitors at the Max Planck Institute for Mathematics in the Sciences, Leipzig;
the Institute of Applied Mathematics, University of Bonn; and the Institute for
Mathematics, University of W\"urzburg.
Some part of this work was performed while TB was visiting the Hausdorff
Research Institute for Mathematics, Bonn.
We acknowledge the hospitality of all these institutions.


\begin{thebibliography}{99}
\bibitem{AC79} Allen, S.\ M.; Cahn, J.\ W. Microscopic theory for
antiphase boundary motion and its application to antiphase domain coarsening,
{\em Acta Metallurgica} 1979, {\bf27}, 1085--1095.
\bibitem{AFP} Ambrosio, L.; Fusco, N.; Pallara, D. {\em Functions of bounded
variation and free discontinuity problems,} Clarendon Press, New York, 2000.
\bibitem{BGN07} Barrett, J.\ W.; Garcke, H.; N\"urnberg, R. A phase field
model for the electromigration of intergranular voids. {\em Interfaces and
Free Boundaries} 2007, {\bf9}, 171--210.
\bibitem{Kaushik} Bhattacharya, K. Comparison of the geometrically nonlinear
and linear theories of martensitic transformation.
{\em Contin.\ Mech.\ Thermodyn.} 1993, {\bf 5}, 205--242.  
\bibitem{BK97} Bhattacharya, K.; Kohn, R.\ V. Elastic Energy Minimization and
the Recoverable Strains of Polycrystalline Shape-Memory Materials.
{\em Arch.\ Rational Mech.\ Anal.} 1997, {\bf139}, 99--180.
\bibitem{BS2} Bhattacharya, K.; Schl\"omerkemper, A.
Stess-induced phase transformations in shape-memory polycrystals,
{\em Arch.\ Rational Mech.\ Anal.} 2010, {\bf 196}, 715--751.
\bibitem{BC11} Blesgen, T.; Chenchiah, I.\ V. A generalized Cahn-Hilliard
equation based on geometrically linear elasticity, 
{\em Interfaces and Free Boundaries} 2011, {\bf 13}, 1--27.
\bibitem{B08} Blesgen, T. The elastic properties of single crystals with
microstructure and applications to diffusion induced segregation.
{\em Crystal Research and Technology} 2008, {\bf43}, 905--913.
\bibitem{BW05} Blesgen, T.; Weikard, U. On the multicomponent Allen-Cahn
equation for elastically stressed solids. {\em Electr.\ J.\ Diff.\ Equ.} 2005,
{\bf89}, 1--17.
\bibitem{BS} Blesgen, T.; Schl\"omerkemper, A.
Towards diffuse interface models with a nonlinear polycrystalline
elastic energy. In: A. Doughett and P. Asnarez (Eds.), Composite
Laminates: Properties, Performance and Applications, Nova Science
Publishing 2010,  pp.\ 465-489.
\bibitem{BCN94} Brochet, D.; Hilhorst, D.; Novick-Cohen, A. Finite-dimensional
exponential attractor for a model for order-disorder and phase separation.
{\em Appl.\ Math.\ Lett.} 1994, {\bf7}, 83--87.
\bibitem{CH} Cahn, J.\ W.; Hilliard, J.\ E.
Free energy of a non-uniform system I. Interfacial free energy.
{\em J. Chem. Phys.} 1958, {\bf28}, 258--267.
\bibitem{CL} Cahn, J.\ W.; Larch{\'e}, F.\ C. The effect of self-stress on
diffusion in solids. {\em Acta Metall.} 1982, {\bf30}, 1835--1845.
\bibitem{CNC94} Cahn, J.\ W.; Novick-Cohen, A. Evolution equations for phase
separation and ordering in binary alloys. {\em J.\ Stat.\ Phys.} 1994,
{\bf 76}, 877--909.
\bibitem{CB2008} Chenchiah, I.\ V.; Bhattacharya, K. The Relaxation of Two-well
Energies with Possibly Unequal Moduli. {\em Arch.\ Rational Mech.\ Anal.} 2008,
{\bf 187}, 409--479.
\bibitem{Che00} Cherkaev, A. V. Variational methods for structural optimization,
{\em Appl. Math. Sci.} 2000, {\bf 140}, Springer publishing.
\bibitem{Donato} Cioranescu, D.; Donato, P. {\em An introduction to
homogenization,} Oxford University Press, 1999.
\bibitem{Dac} Dacorogna, B. {\em Direct methods in the calculus of variations,}
Springer, Berlin, 1989.
\bibitem{Eshelby} Eshelby, J.\ D. Elastic inclusions and inhomogeneities.
{\em Prog.\ Solid Mech.} 1961, {\bf2}, 89--140.
\bibitem{FPL99}  Fratzl, P.; Penrose, O.; Lebowitz, J.\ L. Modelling of phase
separation in alloys with coherent elastic misfit. {\em J.\ Stat.\ Phys.} 1999,
{\bf95}, 1429--1503.
\bibitem{Harald1} Garcke, H. On Cahn-Hilliard systems with elasticity.
{\em Proc.\ Roy.\ Soc.\ Edinburgh Sec.\ A} 2003, {\bf133}, 307--331.
\bibitem{Gottstein} Gottstein, G. Physical Foundations of Materials Science.
Springer Berlin, 2004.
\bibitem{Kohn} Kohn, R.\ V.; Vogelius, M. Relaxation of a variational method
for impedance computed tomography. {\em Comm.\ Pure Appl.\ Math.} 1987,
{\bf60}, 745--777.
\bibitem{Kohn91} Kohn, R.\ V. The relaxation of a double-well energy.
{\em Continuum Mechanics and Thermodynamics} 1991, {\bf3(3)}, 193--236.
\bibitem{KN00} Kohn, R.\ V.; Niethammer, B. Geometrically nonlinear shape-memory
polycrystals made from a two-variant material. {\em Math.\ Mod.\ Num.\ Anal.}
2000, {\bf 34}, 377-398.
\bibitem{Mil01} Milton, G.W. A brief review of the translation method for
bounding effective elastic tensors of composites, {\em Cont. Models and Discrete
Systems} (G.A. Maugin, ed.), Longman Scientific and Technical, 1990, 60--74.
\bibitem{NC00} Novick-Cohen, A. Triple-junction motion for an
Allen-Cahn/Cahn-Hilliard system. {\em Physica D} 2000, {\bf137}, 1--24.
\bibitem{Nye} Nye, J.\ F. {\em Physical properties of crystals: their
representation by tensors and matrices,} Clarendon Press;
New York: Oxford University Press, 1984.
\bibitem{Onuki} Onuki, A. Ginzburg-Landau approach to elastic effects in the
phase separation of solids. {\em J.\ Phys.\ Soc.\ Japan} 1989, {\bf58},
3065--3068.
\bibitem{Salje} Salje, E. Phase transitions in ferroelastic and co-elastic
crystals. Cambridge University Press, 1990.
\bibitem{Ziemer} Ziemer, W. {\em Weakly differentiable functions:
Sobolev Spaces and Functions of Bounded Variation,} Graduate Texts
in Mathematics, Springer, New York, 1989.
\end{thebibliography}
\end{document}